\documentclass[11pt]{article}
\usepackage{epsf}
\usepackage{color}
\usepackage{epsfig}

\usepackage{graphicx}
\usepackage{amsmath}
\usepackage{latexsym}
\usepackage{amssymb}

\usepackage{a4}

\newtheorem{theorem}{Theorem}[section]
\newtheorem{lemma}[theorem]{Lemma}

\numberwithin{equation}{section}

\newcommand{\ad}{{\rm ad}}
\newcommand{\e}{{\rm e}}

\begin{document}

\title{Sufficient conditions for the convergence of the Magnus expansion}

\author{Fernando Casas\thanks{Email: \texttt{Fernando.Casas@uji.es}} \\
   Departament de Matem\`atiques, Universitat Jaume I, \\
  E-12071 Castell\'on, Spain.
   }

\date{}
\maketitle

\begin{abstract}

Two different sufficient conditions are given for the convergence of the Magnus
expansion arising in the study of the linear differential equation $Y' = A(t) Y$. The first one provides
a bound on the convergence domain based on the norm of the operator $A(t)$. The second
condition links the convergence of the expansion with the structure of the spectrum of $Y(t)$, 
thus yielding a more precise characterization. Several examples are proposed to illustrate the main
issues involved and the information on the convergence domain provided by both
conditions.

\vspace*{0.5cm}

\end{abstract}

\section{Introduction}   \label{sec.1}

 The approach followed by Magnus in \cite{magnus54ote}
 to solve the non-autonomous linear differential equation
\begin{equation}   \label{eq.1.1}
    \frac{dY}{dt} = A(t) Y, \qquad Y(0)=I,
\end{equation}
where $Y(t)$ and $A(t)$ are (sufficiently smooth real or complex) $n \times n$ matrices,
is to express the solution $Y(t)$ as the exponential of a certain matrix $\Omega(t)$, 
\begin{equation}   \label{eq.1.2}
   Y(t) = \exp \Omega(t).
\end{equation}
By substituting (\ref{eq.1.2}) into (\ref{eq.1.1}), one can derive the differential equation
satisfied by the exponent $\Omega$ \cite{iserles00lgm}:
\begin{equation}   \label{eq.1.3} 
  \Omega^{\prime} = d \exp_{\Omega}^{-1}(A(t)), \qquad \Omega(0) = O,
\end{equation}
where $d \exp_{\Omega}^{-1}$ is the inverse operator of the power series
\[
    d \exp_{\Omega} = \sum_{k=0}^{\infty} \frac{1}{(k+1)!} \ad_{\Omega}^k \equiv
        \frac{ \exp(\ad_{\Omega}) - I}{\ad_{\Omega}}.
\]
Specifically, its expression is given by    
\[
   d \exp_{\Omega}^{-1}(A) \equiv  \sum_{k=0}^{\infty} \frac{B_k}{k!}
  \mathrm{ad}_{\Omega}^k (A).
\]
Here $\{B_k\}_{k \in \mathbb{Z}_+}$ are the Bernoulli numbers
\cite{abramowitz65hom}, $\mathrm{ad}^k$ is
a shorthand for an iterated commutator,
\[
   \mathrm{ad}_{\Omega}^0 A = A, \qquad \mathrm{ad}_{\Omega}^{k+1} A  =
    [ \Omega, \mathrm{ad}_{\Omega}^k A ],
\]
and  $[\Omega, A] = \Omega A - A \Omega$.  By applying Picard's iteration on
(\ref{eq.1.3}), one gets an infinite series for $\Omega(t)$,
\begin{equation}   \label{eq.1.4}
   \Omega(t) = \sum_{k=1}^{\infty} \Omega_k(t),
\end{equation}     
whose first terms read 
\begin{align}  \label{eq.1.5}
\Omega_{1}(t)  & =\int_{0}^{t}A(t_{1})~\text{d}t_{1},\nonumber\\
\Omega_{2}(t)  & =\frac{1}{2}\int_{0}^{t}\text{d}t_{1}\int_{0}^{t_{1}}%
\text{d}t_{2}\ \left[  A(t_{1}),A(t_{2})\right] \\
\Omega_{3}(t)  & =\frac{1}{6}\int_{0}^{t}\text{d}t_{1}\int_{0}^{t_{1}}%
\text{d}t_{2}\int_{0}^{t_{2}}\text{d}t_{3}\ (\left[  A(t_{1}),\left[
A(t_{2}),A(t_{3})\right]  \right]  +\left[  A(t_{3}),\left[  A(t_{2}%
),A(t_{1})\right]  \right]  )\nonumber
\end{align}
Explicit formulae for $\Omega_m$ of all orders have been given in \cite{iserles99ots},
whereas in \cite{klarsfeld89rgo} a recursive procedure for the generation of any
$\Omega_m$ was proposed, which presents some advantages from a computational
point of view. When this recursion is worked out explicitly, it is possible to express
$\Omega_m$ as a linear combination of $m$-fold integrals of $m-1$ nested commutators
containing $m$ operators $A$, 
\begin{equation}   \label{eq.1.5b}
    \Omega_m(t) =  \sum_{j=1}^{m-1} \frac{B_j}{j!} \,
    \sum_{
            k_1 + \cdots + k_j = m-1 \atop
            k_1 \ge 1, \ldots, k_j \ge 1}
            \, \int_0^t \,
       \ad_{\Omega_{k_1}(s)} \,  \ad_{\Omega_{k_2}(s)} \cdots
          \, \ad_{\Omega_{k_j}(s)} A(s) \, ds    \qquad m \ge 2,
\end{equation}
an expression, however, that becomes increasingly
intricate with $m$, as it should be already evident from the first terms (\ref{eq.1.5}).

Equations (\ref{eq.1.2}) and (\ref{eq.1.4}) constitute the so-called \emph{Magnus
expansion} for the solution of (\ref{eq.1.1}), whereas the infinite series (\ref{eq.1.4}) 
with (\ref{eq.1.5b}) is
known as the \emph{Magnus series}.

Since the 1960s, the Magnus expansion has been successfully applied as a perturbative
tool in numerous areas of physics and chemistry, from atomic and molecular
physics to nuclear magnetic resonance and quantum electrodynamics (see
\cite{blanes98maf}  and \cite{blanes07mem} for a review and a list of references). 
Also, since the work
by Iserles and N{\o}rsett \cite{iserles99ots}, it has been used as a tool to construct
practical algorithms for the numerical integration of equation (\ref{eq.1.1}), while
preserving the main qualitative properties of the exact solution. In this sense, the
corresponding schemes are prototypical examples of \emph{geometric numerical
integrators} \cite{hairer06gni}. 

To be more specific, suppose that $A(t)$ belongs to some matrix Lie 
(sub)\-al\-ge\-bra
$\mathfrak{g}$ for all $t$. Then the exact solution of (\ref{eq.1.1}) evolves
in the  matrix Lie group $\mathcal{G}$ having $\mathfrak{g}$ as its corresponding Lie algebra 
(the tangent space at the identity of $\mathcal{G}$). Observe now that
all terms in the Magnus series are constructed as sums of multiple integrals of nested
commutators, so that  $\Omega$ and indeed any approximation to it obtained by
truncation will also be in the same Lie algebra. Finally, its exponential will be in
$\mathcal{G}$. By truncating appropriately the series, approximating efficiently the
multivariate integrals by quadratures and reducing the number of required commutators,
it is possible to design new integrators based on the Magnus expansion which have
proved to be highly competitive with other, more conventional schemes with respect
to accuracy and computational effort in the numerical integration of (\ref{eq.1.1}) on matrix Lie
groups \cite{blanes00iho,blanes02hoo,iserles00lgm}.

Although the Magnus expansion has been formulated here only for $n \times n$ matrices,
the same result is also valid (as least formally) in a more general setting. As a matter
of fact,  it was originally established assuming only that $A(t)$ is a known
function of $t$ in an associative ring \cite{magnus54ote}. On the other hand, 
in order to apply this approach in Quantum Mechanics,
it is tacitly assumed that the expansion is
also valid when $A(t)$ is a linear operator in a Hilbert space.

From a mathematical point of view, it is clear that there are at least two 
different  issues of paramount importance at the very basis of the Magnus expansion:
\begin{itemize}
 \item First, for what values of $t$ and for what operators $A$ does equation
 (\ref{eq.1.1}) admit a true exponential solution in the form
 (\ref{eq.1.2}) with a certain $\Omega(t)$? This could be called the \emph{existence}
 problem.
 \item Second, given a certain operator $A(t)$, for what values of $t$  does the Magnus
 series (\ref{eq.1.4}) converge? In other words, when $\Omega(t)$ in (\ref{eq.1.2})
 can be obtained as the sum of the series (\ref{eq.1.4})?
 This we describe as the \emph{convergence} problem.
\end{itemize}

Of course, given the relevance of the expansion, both problems have been  
extensively treated in the literature since Magnus proposed this formalism in
1954  \cite{magnus54ote}. 
In section 2 we review some of the most relevant contributions already 
available regarding both aspects, whereas in the rest of the paper we 
will concentrate ourselves on the convergence issue. 
Thus, in section 3 we provide a general
result on the convergence of the Magnus series which is valid for  bounded 
linear operators  $A(t)$  in a Hilbert space, and not only for real matrices. Then, in section 4, 
we analyze the problem from a different point of view,
characterizing the convergence (or divergence) of the series in terms of the 
eigenvalues of the
matrix $Y(t)$. This allows us, in some cases, to obtain more accurate estimates and at
the same time gives us more insight into the convergence problem. Several
examples are also considered to illustrate the main issues involved. 
Finally, section 5 contains a discussion of the results obtained.

\section{Existence and convergence of $\Omega(t)$: previous results}
 \label{sec.2}

\subsection{On the existence of $\Omega(t)$}

In most cases one is interested in the case where 
$A(t)$ belongs to a Lie algebra
$\mathfrak{g}$ under the commutator product. In this rather general
setting, Magnus result can be formulated as four statements
concerning the solution of $dY/dt = A(t) Y$, each one more
stringent than the preceding \cite{wei63nog}. Specifically,
\begin{itemize}
 \item[(A)] The differential equation $dY/dt = A(t) Y$
 has a solution of the form $Y(t) = \exp \Omega(t)$.
 \item[(B)] The exponent $\Omega(t)$ lies in the Lie algebra
  $\mathfrak{g}$.
 \item[(C)] The exponent $\Omega(t)$ is a continuous
 differentiable function of $A(t)$ and $t$, satisfying the
 nonlinear differential equation $d \Omega/dt = d
 \exp_{\Omega}^{-1}(A(t))$.
 \item[(D)] The operator $\Omega(t)$ can be computed by a series
 \[
   \Omega(t) = \Omega_1(t) + \Omega_2(t) + \cdots,
 \]
 where every term is a multivariate integral involving a linear
 combination of nested commutators of $A$ evaluated at different
 times (i.e., the Magnus series (\ref{eq.1.4}) with (\ref{eq.1.5b})).
\end{itemize}
We proceed now to analyze in detail the conditions under which statements
(A)-(C) hold, whereas the validity of (D) will be established by examining the
convergence problem in the rest of the paper.

\vspace*{0.2cm}

\noindent (A) If $A(t)$ and $Y(t)$ are $n \times n$ matrices,
from well-known general theorems on differential equations it is clear that the initial
value problem (\ref{eq.1.1})  always
has a uniquely determined solution $Y(t)$ which is continuous and
has a continuous first derivative in any interval in which $A(t)$
is continuous \cite{coddington55tod}. Furthermore, the determinant
of $Y$ is always different from zero, since
\[
   \det Y(t) = \exp \left( \int_0^t \, \mathrm{tr}\, A(s) ds
   \right).
\]
On the other hand,  any matrix $Y$ can be
written in the form $\exp \Omega$ if and only if $\det Y \ne 0$
\cite[page 239]{gantmacher59tto}, so that it is
always possible to write $Y(t) = \exp \Omega(t)$.

In the general context of Lie groups and Lie algebras, it is
indeed the regularity of the exponential map from the Lie algebra
$\mathfrak{g}$ to the Lie group $\mathcal{G}$ that determines the
global existence of an $\Omega(t) \in \mathfrak{g}$ 
\cite{dixmier57led,saito57scg}: the
exponential map of a complex Lie algebra is globally one to one if
and only if the algebra is nilpotent. In
general, however, the injectivity of the exponential map is only
assured for $\xi \in \mathfrak{g}$ such that $\|\xi\| <
\rho_{\mathcal{G}}$ for a real number $\rho_{\mathcal{G}}
> 0$ and some norm in $\mathfrak{g}$ \cite{moan99ote,moan02obe}. 

\vspace*{0.2cm}

\noindent (B) Although in principle $\rho_{\mathcal{G}}$
constitutes a sharp upper bound for the mere existence of the
operator $\Omega \in \mathfrak{g}$, its practical value in the
case of differential equations is less clear. For instance, 
the logarithm of $Y(t)$ may be complex even for real $A(t)$
 \cite{wei63nog}. In such a situation, the
solution of (\ref{eq.1.1}) cannot be written as the
exponential of a matrix belonging to the Lie algebra over the
field of real numbers. One might argue that this is indeed
possible over the field of complex numbers, but (i) the element
$\Omega$ cannot be computed by the Magnus series (D), since it
contains only real rational coefficients, and (ii) examples exist
where the logarithm of a complex matrix does not lie in the
corresponding Lie subalgebra \cite{wei63nog}. 

It is therefore interesting to determine for which range of $t$ a
real matrix $A(t)$ in (\ref{eq.1.1}) leads to a real logarithm. 
This issue has been tackled in \cite{moan02obe} in the
context of a complete normed (Banach) algebra, proving that if
 \begin{equation}   \label{mo.1}
     \int_0^t \|A(s)\|_2 \, ds < \pi
 \end{equation}
then the solution of (\ref{eq.1.1}) can be written indeed as
 $Y(t) = \exp \Omega(t)$, where $\Omega(t)$ is in the Banach
 algebra. In (\ref{mo.1}), $\| . \|_2$ stands specifically for the 2-norm (or spectral
 norm) of $A$.
 
 \vspace*{0.2cm}

\noindent (C) In his original paper \cite{magnus54ote},
 Magnus was well aware that if
the function $\Omega(t)$ is assumed to be differentiable, it may
not exist everywhere. In fact, he related the differentiability
issue to the existence of the right-hand side of eq. (\ref{eq.1.3}) 
 and gave an implicit condition
for an arbitrary $A$. More specifically, he proved the following
result for the case of $n \times n$ matrices (Theorem V in \cite{magnus54ote}):
\begin{theorem}  \label{con-mag}
  The equation $A(t) = d \exp_{\Omega}(\Omega')$ can be solved by 
  $\Omega' = d \exp_{\Omega}^{-1} A(t)$ for an arbitrary $A$ if and only
  if none of the differences between any two of the eigenvalues of 
  $\Omega$ equals $2\pi i m$, where
 $m= \pm 1, \pm 2,\ldots$, ($m \ne 0$).
\end{theorem}  
  Unfortunately, such a result
has not very much practical application unless we can easily
determine the eigenvalues of $\Omega$ from those of $A(t)$.

\subsection{Convergence of the Magnus series}
\label{convergence}

Let us analyze now in some detail statement (D).
Magnus considered the question of when the series (\ref{eq.1.4})
terminates at some
finite index $m$, thus giving a globally valid $\Omega = \Omega_1
+ \cdots +\Omega_m$. This happens, for instance, if
\[
   \left[ A(t), \int_0^t A(s) ds \right] = 0
\]
identically for all values of $t$, since then $\Omega_k = 0$ for
$k >1$. A sufficient (but not necessary) condition for the
vanishing of all terms $\Omega_k$ with $k > l$ is that
\[
    [A(s_1), [A(s_2), [A(s_3), \cdots ,[A(s_l), A(s_{l+1})] \cdots
    ]]] = 0
\]
for any choice of $s_1, \ldots, s_{l+1}$. In fact, the termination
of the series cannot be established solely by consideration of the
commutativity of $A(t)$ with itself, and Magnus considered an
example illustrating this point \cite{magnus54ote}.

In general, however, the Magnus series does not converge unless
$A$ is small in a suitable sense. Several bounds to the actual
radius of convergence in terms of $A$ have been obtained along the years. 
Most of these results can be stated as follows. 
If $\Omega_m(t)$ denotes the homogeneous element with $m$
commutators in the Magnus series as given by (\ref{eq.1.5b}),
then $\Omega(t) = \sum_{m=1}^{\infty} \Omega_m(t)$ is absolutely
convergent for $0 \le t < T$, with
\begin{equation}   \label{eq.2.1}
    T = \max \left\{ t \ge 0 \, : \, \int_0^t \|A(s)\|_2 \, ds < r_c
      \right\}.
\end{equation}
Thus, Pechukas and Light \cite{pechukas66ote} and
Karasev and Mosolova \cite{karasev77ipa} both obtained $r_c=\log 2 =
0.693147\ldots$, whereas Chacon and Fomenko \cite{Chacon91rff} got
$r_c=0.57745\ldots$. In 1998, Blanes \textit{et al.}
\cite{blanes98maf} and Moan \cite{moan98eao} obtained independently
the improved bound
\[
    r_c = \frac{1}{2}  \int_0^{2 \pi} \frac{1}{2 + \frac{x}{2} (1 - \cot \frac{x}{2})}
    \,  dx  \equiv \xi = 1.08686869\ldots
\]
Based on the analysis of some selected examples, Moan \cite{moan02obe}
concluded that, in order to get convergence \emph{for all real matrices}
$A(t)$, necessarily $r_c \le \pi$, and more recently Moan and Niesen
\cite{moan06cot} have been able to prove rigorously that indeed $r_c=\pi$. 

This result shows, in particular, that statement (D) is locally valid,
but cannot be used to compute $\Omega$ in the large. However, as
we have seen, the other statements need not depend on the validity
of (D). In particular, if (B) and (C) are globally valid, one can
still investigate many of the properties of $\Omega$ even though
one cannot compute it with the aid of (D).

\section{A generic result on the convergence of the Magnus series}
 \label{sec.3}

 \subsection{General formulation}
 
 As we have mentioned before, if $A(t)$ is a real $n \times n$ matrix, then
(\ref{mo.1}) gives a condition for $Y(t)$ to have a real logarithm. Moreover, it
has been shown that,
under the same condition, the Magnus series (\ref{eq.1.4}) converges precisely
to this logarithm, i.e., its sum $\Omega(t)$ satisfies $\e^{\Omega(t)} = Y(t)$
\cite{moan06cot}. Our purpose in this section is provide a different proof 
of this property which in fact is
also valid in the more general setting of linear operators in a Hilbert space
of arbitrary dimension.

To begin with, let $A(t)$ be a bounded linear operator in a Hilbert space $\mathcal{H}$,
with $2 \le \mathrm{dim} \ \mathcal{H} \le \infty$. Let us introduce a new parameter
$\varepsilon \in \mathbb{C}$ and denote by $Y(t;\varepsilon)$  the 
solution of the initial value problem
\begin{equation}   \label{eq.3.1}
   \frac{dY}{dt} = \varepsilon A(t) Y,  \qquad Y(0)=I,
\end{equation}
where now $I$ denotes the identity operator in $\mathcal{H}$.  It is known that  
$Y(t;\varepsilon)$ is an analytic function of $\varepsilon$ for a fixed value of $t$. 
Let us introduce the set $B_{\gamma} \subset \mathbb{C}$ 
characterized by the real parameter $\gamma$,
\[
  B_{\gamma} = \{ \varepsilon \in \mathbb{C} \, : \, |\varepsilon| 
         \int_0^t \|A(s)\| \, ds <  \gamma \}.
\]            
Here $\|.\|$ stands for the norm defined by the inner product on $\mathcal{H}$.    
Our first statement is that, if $t$ is fixed, the operator function 
$\varphi(\varepsilon) = \log Y(t;\varepsilon)$ is well defined
in $B_{\gamma}$ when $\gamma$ is small enough, say $\gamma < \log 2$, as an
analytic function of $\varepsilon$. 

As a matter of fact, this is a direct consequence of the results collected in section \ref{convergence}:
if, in particular, $|\varepsilon| \int_0^t \|A(s)\| \, ds < \log 2$, the 
Magnus series corresponding to (\ref{eq.3.1}) converges and its sum
$\Omega(t; \varepsilon)$ satisfies $\e^{\Omega(t; \varepsilon)} = Y(t; \varepsilon)$.
In other words, the power series $\Omega(t; \varepsilon)$ coincides with
$\varphi(\varepsilon)$ when $|\varepsilon| \int_0^t \|A(s)\| \, ds < \log 2$, and so
the Magnus series is the power series expansion of $\varphi(\varepsilon)$ around
$\varepsilon = 0$.

The next theorem shows that, indeed, $\gamma = \pi$.
\begin{theorem}   \label{main-theorem}
       The function $\varphi(\varepsilon) = \log Y(t;\varepsilon)$ is an analytic 
       function of $\varepsilon$ in the set $B_{\pi}$, with
\[
  B_{\pi} = \{ \varepsilon \in \mathbb{C} \, : \, |\varepsilon| 
         \int_0^t \|A(s)\| ds <  \pi \}.
\]                       
If $\mathcal{H}$ is infinite-dimensional, the statement holds true if $Y$ is a normal
operator.
\end{theorem}

The proof of this theorem is based on some elementary properties of the unit
sphere $S^1$ in a Hilbert space. Let us define the angle between any two vectors
$x \ne 0$, $y \ne 0$ in $\mathcal{H}$, $\mathrm{Ang}\{x,y\} = \alpha$, $0 \le \alpha \le \pi$, from
\[
      \cos \alpha = \frac{\mathrm{Re} \langle x, y \rangle}{\|x\| \, \|y\|},
\]    
where $\langle \cdot, \cdot \rangle$ is the inner product on $\mathcal{H}$.
This angle is a metric in $S^1$, i.e., the triangle inequality holds there. 
A trivial property which will be used in the sequel is the following: if
$\|x\| = 1$ and $\|u\| \le 1/2$, then $\mathrm{Ang}\{x+u,x\} \le  \|u\| (1 + \|u\|^2)$.

The second basic property of the angle we need
is given by the following lemma, whose proof (due to Moan \cite{moan02obe}) is included here for
completeness.
\begin{lemma}  \label{lemma1}
 (Moan). For any $x \ne 0$ in $\mathcal{H}$, 
$\mathrm{Ang}\{Y(t;\varepsilon) x, x\} \le |\varepsilon| \int_0^t \|A(s)\| ds$
\end{lemma}  
\noindent \textit{Proof of Lemma \ref{lemma1}}. Let $y_0 \equiv x$ and consider
the vector $y(t) = Y(t;\varepsilon) y_0$ satisfying the initial value problem
$y^{\prime} = \varepsilon A(t) y$, $y(0)=y_0$. Then, clearly,
$\|y^{\prime}\| \le |\varepsilon| \|A(t)\| \, \|y\|$. Let $\hat{y}(t) = \frac{y(t)}{\|y(t)\|}$ denote
the unit vector in the direction of $y(t)$, so that $y^{\prime} = \frac{d \|y\|}{dt}  \, \hat{y} + 
 \|y\| \hat{y}^{\prime}$. On the other hand, since $\langle \hat{y}, \hat{y} \rangle = 1$, then
 $\langle \hat{y}^{\prime}, \hat{y} \rangle + \langle \hat{y}, \hat{y}^{\prime} \rangle =
 \mathrm{Re} \langle \hat{y}, \hat{y}^{\prime} \rangle = 0$, i.e., $\hat{y}$ and 
 $\hat{y}^{\prime}$ are orthogonal. In consequence,
 \[
   \langle y^{\prime}, y^{\prime} \rangle = \|y^{\prime}\|^2 = \left(
     \frac{d \|y\| }{dt} \right)^2 + \|y\|^2  \, \|\hat{y}^{\prime}\|^2,
 \]    
whence, by discarding the $(\|y\|^\prime)^2$ term,   
$\|\hat{y}^{\prime}\|^2 \, \|y\|^2 \le \|y^{\prime}\|^2$, and thus
$\|\hat{y}^{\prime}\| \, \|y\| \le \|y^{\prime}\| \le |\varepsilon| \|A(t)\| \, \|y\|$ or
simply $\|\hat{y}^{\prime}\| \le |\varepsilon| \|A(t)\|$. Integrating this last inequality we get
\[
   \int_0^t \|\hat{y}^{\prime}(s)\| \, ds \le |\varepsilon| \int_0^t \|A(s)\| ds,
\]
but
\[
       \int_0^t \|\hat{y}^{\prime} (s)\| \, ds = \int_0^t \sqrt{ \langle \hat{y}'(s), \hat{y}'(s)  \rangle} \, ds
\]
is the length (defined through the metric given by the inner product) of the curve traced       
by the unit vector $\hat{y}(s)$ when $s \in [0,t]$ on the unit sphere $S^1$, which is 
greater than or equal to
$\mathrm{Ang}\{y(t),y_0\}$, and this proves the result.
\hfill{$\Box$}

Observe that if $Y$ is a normal operator in $\mathcal{H}$, i.e., $Y Y^* = Y^* Y$, where
$Y^*$ denotes the adjoint operator of $Y$ (in particular, if $Y$ is unitary), then
$\|Y^* x\| = \|Yx\|$ for all $x \in \mathcal{H}$ and therefore 
$\mathrm{Ang}\{Y^* x, x\} = \mathrm{Ang}\{Yx, x\}$.

The following lemma provides useful information on the location of the 
eigenvalues of a given bounded linear operator in $\mathcal{H}$ \cite{mityagin90unp}.
  \begin{lemma}  \label{lemma2}
    (Mityagin). Let $T$ be a (bounded) operator on $\mathcal{H}$. 
    If $\mathrm{Ang}\{T x,x\} \le \gamma$ and $\mathrm{Ang}\{T^* x,x\} \le \gamma$
    for any $x \ne 0$, $x \in \mathcal{H}$, where $T^*$ denotes the adjoint operator of $T$,then the spectrum of $T$, $\sigma(T)$, is contained in the set
    \[
       \Delta_{\gamma} = \{ z = |z| \e^{i \omega} \in \mathbb{C} \, : \, |\omega| \le \gamma\}
    \]   
  \end{lemma}
\noindent \textit{Proof of Lemma \ref{lemma2}}. Without loss of generality, we may assume
$\gamma < \pi$ (if $\gamma \ge \pi$, there is no statement here). If
$\mathrm{dim} \ \mathcal{H} < \infty$, only the first requirement on $T$, 
$\mathrm{Ang}\{T x,x\} \le \gamma$ is sufficient, since in that case,
if $\lambda =
|\lambda| \, \e^{i \omega} \ne 0$, $-\pi < \omega \le \pi$, is in $\sigma(T)$, then there exists some
$f \ne 0$ such that $T f = \lambda f$ and 
\[
   \mathrm{Ang}\{T f,f\} = \mathrm{Ang}\{\lambda f, f\} = |\omega| \le \gamma.
\]   
If, on the other hand, $\mathrm{dim} \ \mathcal{H} = \infty$ and $\lambda \in \sigma(T)$,
$\lambda \ne 0$, then, as is well known, either (i) $\lambda$ belongs to the approximate
spectrum of $T$, $\sigma_{ap}(T)$, or (ii) $\lambda$ is in the residual spectrum,
$\sigma_r(T)$ \cite{hunter01aan}.

\noindent (i) In the first case, there is a sequence $\{f_n\}$ in $\mathcal{H}$ such that
$\|f_n\| = 1$ for all $n$ and $\lim_{n \rightarrow \infty} \|(T-\lambda I) f_n \| = 0$. Equivalently,
$Tf_n = \lambda f_n + \varepsilon_n$, with $\|\varepsilon_n\| \rightarrow 0$ when
$n \rightarrow \infty$. Then we have
\begin{equation}  \label{dem.1}
  \gamma \ge \mathrm{Ang}\{T f_n, f_n\} = \mathrm{Ang}\{\lambda f_n + \varepsilon_n, f_n\}
    \ge \mathrm{Ang}\{\lambda f_n, f_n\} - \mathrm{Ang}\{\lambda f_n, \lambda f_n +
     \varepsilon_n\}
\end{equation}
since the angle is a metric in $S^1$. Now, as 
$\lambda = |\lambda| \, \e^{i \omega} \ne 0$, $-\pi < \omega \le \pi$, it is clear from
(\ref{dem.1}) that
\begin{equation}   \label{dem.2}
   \gamma \ge |\omega| - \mathrm{Ang}\{f_n, f_n + \frac{1}{\lambda} \varepsilon_n\} \ge
      |\omega| - \frac{\|\varepsilon_n\|}{|\lambda|} \left( 1+ \frac{\|\varepsilon_n\|^2}{|\lambda|^2}
      \right),
\end{equation}      
where the last inequality holds when  $\frac{\|\varepsilon_n\|}{|\lambda|} \le \frac{1}{2}$, i.e.,
for sufficiently large $n$. Taking the limit $n \rightarrow \infty$ in (\ref{dem.2}) leads
to $\gamma \ge |\omega|$.

\noindent (ii) If $\lambda \in \sigma_r(T)$, then $\bar{\lambda}$ is an eigenvalue of $T^*$,
i.e., $\mathrm{Ker} \, (T^* - \bar{\lambda} I) \ne \{0\}$ \cite{hunter01aan}. Since, by
assumption, $\mathrm{Ang}\{T^* x,x\} \le \gamma$ for all $x \ne 0$, we can apply again the
argument in (i) to $T^*$, $\bar{\lambda}$ and conclude that $|\omega| \le \gamma$.
\hfill{$\Box$}

\noindent Now we are ready to prove the main theorem.

\noindent \textit{Proof of Theorem \ref{main-theorem}}. Let us introduce the 
operator $T \equiv Y(t;\epsilon)$, with $\varepsilon \in B_{\gamma}$, $\gamma < \pi$. Then
by Lemma \ref{lemma1}, $\mathrm{Ang}\{T x,x\} \le \gamma$ for all $x \ne 0$, and thus, by
Lemma \ref{lemma2}, 
\begin{equation}   \label{eq.pr.1}
    \sigma(T) \subset  \Delta_{\gamma}.
\end{equation}    
If $\mathrm{dim} \ \mathcal{H} = \infty$ and we assume that $Y(t;\epsilon)$ is a normal operator,
then (\ref{eq.pr.1}) also holds.

From equation (\ref{eq.3.1}) in integral form,
\[
  Y(t; \varepsilon) = I + \varepsilon \int_0^t A(s) Y ds,
\]
one gets $\|Y\| \le 1 + |\varepsilon| \int_0^t \|A(s)\| \, \|Y\| ds$, and application of Gronwall's
lemma \cite{gronwall19not} leads to
\[
    \|Y(t;\varepsilon)\| \le \exp \left( |\varepsilon| \int_0^t \|A(s)\| ds \right).
\]
An analogous reasoning for the inverse operator also proves that      
\[
    \|Y^{-1}(t;\varepsilon)\| \le \exp \left( |\varepsilon| \int_0^t \|A(s)\| ds \right).
\]
In consequence,
\[
   \|T\| \le \e^{\gamma}  \qquad \mbox{ and } \qquad \|T^{-1}\| \le \e^{\gamma}.
\]

If $\lambda \ne 0 \in \sigma(T)$, then $|\lambda| \le \|T\|$ \cite{hunter01aan}
and therefore $|\lambda| \le \e^{\gamma}$. In addition, $\frac{1}{\lambda} \in
\sigma(T^{-1})$, so that $|\lambda| \ge \e^{-\gamma}$. Equivalently,
\begin{equation}  \label{eq.pr.2}
  \sigma(T) \subset \{ z \in \mathbb{C} : \e^{-\gamma} \le |z| \le \e^{\gamma} \} \equiv G_{\gamma}.   
\end{equation}
Putting together (\ref{eq.pr.1}) and (\ref{eq.pr.2}), one has
\[
   \sigma(T) \subset G_{\gamma} \cap \Delta_{\gamma} \equiv \Lambda_{\gamma}.
\]
Now choose any value $\gamma'$ such that  $\gamma < \gamma' < \pi$ 
(e.g., $\gamma' = (\gamma + \pi)/2$) and consider the closed curve 
$\Gamma = \partial \Lambda_{\gamma'}$.  Notice that the curve $\Gamma$
encloses $\sigma(T)$ in its interior, so that it is possible to define the function
$\varphi(\varepsilon) = \log Y(t; \varepsilon)$ by the equation 
\cite{dunford58lop}
 \begin{equation}   \label{eq.pr.3}
      \varphi(\epsilon) = \frac{1}{2 \pi i} \int_{\Gamma} \log z \, (z I - Y(t;\epsilon))^{-1}  \, dz,
  \end{equation}
where the integration along $\Gamma$ is performed in the counterclockwise direction.
As is well known, (\ref{eq.pr.3}) defines an analytic function of $\varepsilon$ in
$B_{\gamma'}$ \cite{dunford58lop} and the result of the theorem follows.
\hfill{$\Box$}

\begin{theorem}  \label{conv-mag}
  Let us consider the differential equation $Y' = A(t) Y$ defined in 
  a Hilbert space $\mathcal{H}$ with $Y(0)=I$, and let $A(t)$ be a bounded linear operator
  on $\mathcal{H}$. Then, the Magnus series 
  $\Omega(t) = \sum_{k=1}^{\infty}  \Omega_k(t)$, with
$\Omega_k$ given by (\ref{eq.1.5b}) converges  in the interval $ t \in [0,T)$  such that
\[
   \int_0^T \|A(s)\| ds < \pi
\]
and the sum $\Omega(t)$ satisfies $\exp \Omega(t) = Y(t)$. The statement also holds when
$\mathcal{H}$ is infinite-dimensional if $Y$ is a normal operator (in particular, if $Y$ is
unitary).
\end{theorem}    

\noindent \textit{Proof}.
Theorem \ref{main-theorem} shows that $\log Y(t;\varepsilon) \equiv \varphi(\varepsilon)$ is a 
 well defined and analytic function of $\varepsilon$ for 
  \[
       |\varepsilon| \int_0^t \|A(s)\| ds  < \pi.
  \]     
It has also been shown that the Magnus series 
$\Omega(t;\varepsilon) = \sum_{k=1}^{\infty} \varepsilon^k \Omega_k(t)$, with
$\Omega_k$ given by (\ref{eq.1.5b}), is absolutely convergent when
$|\varepsilon| \int_0^t \|A(s)\| ds  < \xi = 1.0868...$ and its sum satisfies
$\exp \Omega(t;\varepsilon) = Y(t;\varepsilon)$. Hence, the Magnus series is
the power series of the analytic function $\varphi(\varepsilon)$ in the disk 
$|\varepsilon| < \xi / \int_0^t \|A(s)\| ds$. But $\varphi(\varepsilon)$ is analytic
in $B_{\pi} \supset B_{\xi}$ and the power series has to be unique. In consequence,
the power series of $\varphi(\varepsilon)$ in $B_{\pi}$ has to be same as the power
series of $\varphi(\varepsilon)$ in $B_{\xi}$, which is precisely the Magnus series. Finally,
by taking $\varepsilon = 1$ we get the desired result.
\hfill{$\Box$}

\subsection{Examples}
\label{examples1}
 
Theorem \ref{conv-mag} provides thus sufficient conditions for the convergence 
of the Magnus series based on an estimate by the norm of the operator $A$. In
particular, it guarantees that the operator $\Omega(t)$  in $Y(t) = \exp \Omega(t)$
can safely be obtained with the convergent series $\sum_{k \ge 1} \Omega_k(t)$ for
$0 \le t < T$ when the terms $\Omega_k(t)$ are computed with (\ref{eq.1.5b}).   A natural
question arising here is the following: is the bound estimate provided by Theorem
\ref{conv-mag} sharp or is there still room for improvement? In order to clarify this issue,
we next analyze two simple examples involving $2 \times 2$ matrices.

\vspace*{0.3cm}

\noindent \textbf{Example 1}. Moan and Niesen \cite{moan06cot} consider the 
initial value problem (\ref{eq.1.1}) with
 \begin{equation}   \label{ej1.1}
   A(t) = \left(  \begin{array}{rr}
        2  &  t \\
        0  & -1
           \end{array}   \right).
 \end{equation}
If we introduce, as before, the complex parameter $\varepsilon$ in the problem, the
corresponding exact solution $Y(t;\varepsilon)$ of (\ref{eq.3.1}) is given by
 \begin{equation}   \label{ej1.2}
   Y(t;\varepsilon) = \left(  \begin{array}{lc}
        \e^{2 \varepsilon t}   &  \ \ \frac{1}{9 \varepsilon} \e^{2 \varepsilon t} - \left(
           \frac{1}{9 \varepsilon} + \frac{1}{3} t \right) \e^{-\varepsilon t}   \\
        0  &   \e^{-\varepsilon t} 
           \end{array}   \right)
 \end{equation}
and therefore
   \[
       \log Y(t;\varepsilon) = \left(  \begin{array}{rc}
        2t  &  g(t;\varepsilon) \\
        0  & -t
           \end{array}   \right), \quad \mbox{ with } \ \ g(t;\varepsilon) = \frac{t (1 - \e^{3 \varepsilon t} + 
           3 \varepsilon t)}{3(1- \e^{3 \varepsilon t})}.
    \]           
The Magnus series can be obtained by computing the Taylor expansion of  $\log Y(t;\varepsilon)$ 
around $\varepsilon = 0$. Notice that the function $g$ has a singularity when
$\varepsilon t = \frac{2\pi}{3} i$, and thus, by taking  $\varepsilon = 1$, the Magnus series only converges up to $t=  \frac{2}{3} \pi$. On the other hand, condition  
$\int_0^T \|A(s)\| ds  < \pi$ leads to $T \approx 1.43205 < \frac{2}{3} \pi$.
In consequence, the actual convergence domain of the Magnus series
is larger than the estimate provided by Theorem \ref{conv-mag}.

\vspace*{0.3cm}

\noindent \textbf{Example 2}. Let us introduce the matrices
\begin{equation}   \label{ej2.1}
    X_{1}=  \left(
  \begin{array}{cr}
  1  & 0 \\
  0  &  -1
  \end{array} \right), \qquad
  X_{2}=  \left(
  \begin{array}{cc}
  0  & 1 \\
  0  & 0
  \end{array} \right)
\end{equation}
and define 
\[
  A(t)= \left\{
  \begin{array}{ll}
  \beta \, X_{2} \quad & 0\leq t \leq 1 \\
  \alpha \, X_{1} \quad & t > 1
  \end{array} \right.
\]  
with $\alpha, \beta$ complex constants. Then, the solution of
equation (\ref{eq.1.1}) at $t=2$ is $Y(2) = \e^{\alpha X_{1}} \e^{\beta X_{2}}$,
so that
\begin{equation}   \label{ex1.4}
 \Omega(2) = \log(\e^{\alpha X_{1}} \e^{\beta X_{2}})  = \alpha X_{1} +
    \frac{2\alpha \beta}{1- \e^{-2\alpha}} \, X_{2},
\end{equation}
an analytic function if $|\alpha| < \pi$ with first singularities
at $\alpha = \pm i \pi$.

On the other hand, a simple calculation with the recurrence (\ref{eq.1.5b}) shows that
\begin{equation}  \label{ex2.2}
  \Omega(2) = \sum_{k=1}^{\infty} \Omega_{k}(2) =
     \alpha X_{1} + \beta X_{2} + \sum_{n=2}^{\infty}
    (-1)^{n-1}\frac{2^{n-1}B_{n-1}}{(n-1)!} \ \alpha^{n-1}
     \beta \, X_{2}.
\end{equation}
Comparing with expression (\ref{ex1.4}), it is clear that the Magnus series
 cannot converge at $t=2$ if $|\alpha| \ge \pi$, independently
of $\beta \ne 0$. 

If we take the spectral norm, then $\|X_1\| = \|X_2\| = 1$ and
\[
   \int_0^{t=2} \|A(s)\| ds = |\alpha| + |\beta|,
\]
so that the convergence domain provided by
Theorem \ref{conv-mag} is $|\alpha| + |\beta| < \pi$ for this example.

\vspace*{0.4cm}

From the analysis of Examples 1 and 2 we can conclude the following. First,
the convergence domain of the Magnus series provided by Theorem \ref{conv-mag}
is the best result one can get for a generic bounded operator $A(t)$ in a Hilbert
space, in the sense that one may consider specific $A(t)$, as in Example 2, where
the series diverges for any time $t$ such that  $\int_0^t \|A(s)\| ds > \pi$. Second, there
are also situations (as in Example 1) where the bound estimate $r_c=\pi$ 
is still rather conservative: \emph{the  Magnus series converges indeed for a 
larger time interval than that
given by Theorem \ref{conv-mag}}. This is particularly evident if one considers equation
(\ref{eq.3.1}) with $A(t)$ a diagonal matrix, 
 \begin{equation}   \label{diag1}
   A(t) = \left(  \begin{array}{cc}
        a_1(t)  &  0 \\
        0  &  a_2(t)
           \end{array}   \right).
 \end{equation}
Then, the exact solution $Y(t;\varepsilon)$ of (\ref{eq.3.1}) is a diagonal matrix whose elements
are non-vanishing entire functions of $\varepsilon$, and obviously $\log Y(t;\varepsilon)$
is also an entire function of $\varepsilon$. In such circumstances, the convergence domain
 $|\varepsilon| \int_0^t \|A(s)\| ds  < \pi$ for the Magnus series does not make much sense.

\section{Another characterization of the convergence of the Magnus series}

\subsection{Main result on convergence}

The examples collected in the preceding section (and many more one can build) clearly
show that, although the condition $\int_0^t \|A(s)\| ds  < \pi$ is sharp (in the sense that
the constant $\pi$ is the largest number for which Theorem \ref{conv-mag}  holds in general), it is certainly  not necessary for the convergence of the Magnus series. Thus, it would
be highly desirable to have a more realistic criterion  which give both necessary 
and sufficient conditions for convergence.

In \cite{moan06cot}, a conjecture is formulated, relating the convergence of the Magnus
series with the eigenvalues of the exact solution $Y(t;\varepsilon)$. Here 
we state a theorem which, on the one hand, explains the phenomena observed by
Moan and Niesen \cite{moan06cot} and, on the other hand, provides a new 
tool to determine the actual convergence domain of the Magnus series in some physically
relevant examples and applications.

The main result in this section (Theorem \ref{thY-S}) 
is valid for complex $n \times n$ matrices and is based on the theory
of analytic matrix functions, in particular, in the logarithm of an analytic matrix function.
In fact, it is a direct consequence of the analysis done in
 \cite[Chapter 1, section 3]{yakubovich75lde}. Here we shall summarize the most
 relevant aspects of the formalism and refer the reader to \cite{yakubovich75lde} for 
 a more detailed treatment (including proofs).

Our starting point is again the initial value problem $Y' = \varepsilon A(t) Y$,
$Y(0)=I$, where now $A(t)$ and $Y$ are (complex) $n \times n$ matrices and
$\varepsilon \in \mathbb{C}$. If we denote by $Y_t(\varepsilon)$ the exact solution for
a fixed value of $t$, $Y_t(\varepsilon) \equiv Y(t;\varepsilon)$, it is clear that 
$Y_t(\varepsilon)$ is an analytic function of $\varepsilon$ \cite{coddington55tod}, since
the Neumann series
\[
    Y_t(\varepsilon) = I + \sum_{k=1}^{\infty} \varepsilon^k 
    \int_0^t dt_1 A(t_1) \int_0^{t_{1}} dt_2 A(t_2) \cdots \int_0^{t_{k-1}} dt_k A(t_k)      
\]
converges provided that $\int_0^t \|A(s)\| ds < \infty$. In addition, 
$\det Y_t(\varepsilon) \ne 0$ for all $\varepsilon$. Under these conditions, it has been shown
that the matrix $\Omega_t(\varepsilon) = \log Y_t(\varepsilon)$ is also an analytic function
of $\varepsilon$ at $\varepsilon = 0$. In other words, the series
$\Omega_t(\varepsilon) = \sum_{k \ge 1} \varepsilon^k \Omega_{t,k}$ (i.e., the Magnus
series) is convergent for
sufficiently small  $\varepsilon$. The goal is then to determine the actual radius 
of convergence $r$ of this series. 

Let us denote by $\rho_1(\varepsilon), \ldots, \rho_n(\varepsilon)$ the eigenvalues
  of the matrix $Y_t(\varepsilon)$. Notice that $Y_t(0) = I$, so that $\rho_1(0) = \cdots = \rho_n(0) = 1$.
It is therefore natural to take the principal values of the logarithm,
$ \log \rho_1(0) = \cdots = \log \rho_n(0) = 0$, as this choice is consistent with the series
$\Omega_t(\varepsilon)$. 

Let $L$ be a curve on the $\varepsilon$ plane in the disk $|\varepsilon| < r_0 < \infty$ issuing
from the origin. Recall that the matrix $Y_t(\varepsilon)$ is analytic in the disk 
$|\varepsilon| < r_0$. On the curve $L$ it is possible to define a unique function
$\log \rho_j(\varepsilon)$, $j=1,\ldots,n$, by continuity, given the values
$\log \rho_j(0) = 0$.  

Let $\rho_0$ be a multiple eigenvalue of $Y_t(\varepsilon_0)$
for some $\varepsilon_0$ with $|\varepsilon_0| < r_0$ with multiplicity $l$. 
If we reorder the eigenvalues of
$Y_t(\varepsilon_0)$ in such a way that the first $l$ are precisely $\rho_0$, it is clear that
 the numbers $\log \rho_1(\varepsilon_0)$, $\log \rho_2(\varepsilon_0)$, $\ldots, 
\log \rho_{l}(\varepsilon_0)$, $1 < l \le n$, are
 congruent modulo $2 \pi i$ and are such that 
 $\rho_1(\varepsilon_0) = \cdots = \rho_{l}(\varepsilon_0) = \rho_0$. Associated with
 this multiple eigenvalue $\rho_0$ there is a pair of integers $(p,q)$ defined as
 follows.

The integer $p$ is the greatest number of equal terms in the set of numbers 
$\log \rho_1(\varepsilon_0)$, $\log \rho_2(\varepsilon_0)$, $\ldots, 
\log \rho_{l}(\varepsilon_0)$ such that  $\rho_k(\varepsilon_0) = \rho_0$, $k=1,\ldots, l$. 

The integer $q$ is the maximum degree of the elementary divisors $(\rho - \rho_0)^k$ of 
$Y_t(\varepsilon_0)$,
i.e., the maximum dimension of the elementary Jordan block corresponding to $\rho_0$.

Notice that the numbers $l$ and $q$ depend only on the particular eigenvalue $\rho_0$,
whereas the integer $p$ depends on $\rho_0$ and the curve $L$.

Under these conditions, it is possible to prove the following lemma \cite[page 64]{yakubovich75lde}
on the convergence of the series $\Omega_t(\varepsilon)$.
\begin{lemma}  \label{yaku}
 (Yakubovich--Starzhinskii). Suppose that the series \newline
 $\Omega_t(\varepsilon) = \sum_{k \ge 1} \varepsilon^k \Omega_{t,k}$ satisfies  that
$\exp \Omega_t(\varepsilon) = Y_t(\varepsilon)$ for sufficiently small $|\varepsilon|$. Then
\begin{itemize}
 \item[(a)] If $r < r_0$ is the radius of convergence of the series
$\Omega_t(\varepsilon)$, the eigenvalues $\lambda_1(\varepsilon), \ldots,
\lambda_n(\varepsilon)$ of the matrix $\Omega_t(\varepsilon)$, defined for $|\varepsilon| < r$,
can be defined by continuity on the circle $|\varepsilon| = r$, and there exists a point
$\varepsilon_0$ such that for some $j,k=1,\ldots, n$
\[
    \lambda_j(\varepsilon_0) - \lambda_k(\varepsilon_0) = 2 \pi i m,
\]
where $m \ne 0$ is an integer.
 \item[(b)] Suppose that $\varepsilon_0$ is the value of $\varepsilon$ of smallest absolute
 value ($\varepsilon_0 \ne 0$, $|\varepsilon_0| < r_0$) such that the matrix $Y_t(\varepsilon_0)$
 has an igenvalue $\rho_0$ of multiplicity $l > 1$. Suppose that there is at least one
 such an eigenvalue $\rho_0$ and at least one curve in the disk $|\varepsilon| < |\varepsilon_0|$
 joining the origin $\varepsilon = 0$ with the point $\varepsilon = \varepsilon_0$ such that
 $p < q$, where the integers $p$ and $q$ have been defined before. Then $r = |\varepsilon_0|$
 is the radius of convergence of the series 
 $\Omega_t(\varepsilon) = \sum_{k \ge 1} \varepsilon^k \Omega_{t,k}$.
\end{itemize}
\end{lemma}

In order to apply this result one first has to solve the equation
\begin{equation}  \label{eq.4.1}
   \Delta(\varepsilon) = 0,
\end{equation}   
where $\Delta(\varepsilon)$ denotes the discriminant of the characteristic polynomial
$\det(Y_t(\varepsilon) - \rho I)$. We recall here that the discriminant of a polynomial 
\[
    p(x) = a_n x^n + a_{n-1} x^{n-1} + \cdots + a_1 x + a_0
 \]
is given by
 \[
    a_n^{2n-2} \prod_{i < j} (r_i - r_j)^2,
 \]
 with $r_1, \ldots, r_n$ complex roots of $p(x)$, so that it vanishes if and only if $p(x)$ has
 one or more multiple roots \cite{kurosh72hal}. Thus, it can be used to test for the presence
 of multiple roots, without having to actually compute the roots of $p(x)$.

We write the solutions of equation (\ref{eq.4.1}) in order of non-decreasing absolute
value,
\begin{equation}  \label{eq.4.2}
   \varepsilon_0^{(1)}, \, \varepsilon_0^{(2)}, \varepsilon_0^{(3)}, \ldots
\end{equation}
and consider the circle $|\varepsilon| = |\varepsilon_0^{(1)}|$ in the complex 
$\varepsilon$-plane. Let $\rho_0^{(1)}$ denote an  eigenvalue of 
$Y_t(\varepsilon_0^{(1)})$ with multiplicity $l_1 > 1$. Let $\varepsilon$ move along some fixed curve 
$L$ from $\varepsilon = 0$ to 
$\varepsilon = \varepsilon_0^{(1)}$ in the circle $|\varepsilon| \le |\varepsilon_0^{(1)}|$. Then
it is clear that  $l_1$ eigenvalues $\rho_j(\varepsilon)$ will tend to $\rho_0^{(1)}$ at 
$\varepsilon = \varepsilon_0^{(1)}$. If these points lie at 
$\varepsilon = \varepsilon_0^{(1)}$  on the same sheet of the Riemann
surface of the function $\log z$, and this is true for all (possible) multiple eigenvalues 
of $Y_t(\varepsilon)$ at $\varepsilon = \varepsilon_0^{(1)}$, then
$\varepsilon_0^{(1)}$ is called a \emph{extraneous root} of equation
(\ref{eq.4.1}). Otherwise, 
$\varepsilon_0^{(1)}$ is called a \emph{non-extraneous root}.

Now, by Lemma \ref{yaku}, when $|\varepsilon| < |\varepsilon_0^{(1)}|$, the series
for $\Omega_t(\varepsilon)$ is convergent, 
so that the numbers $\log \rho_j(\varepsilon)$ are uniquely determined up to
multiplicity as eigenvalues of the matrix $\Omega_t(\varepsilon)$.

If $\varepsilon_0^{(1)}$ is an extraneous root, there is no obstacle to the 
convergence of the series and thus we proceed to the next value in the
sequence (\ref{eq.4.2}). We continue this classification until
 a non-extraneous root is obtained. Assume, for simplicity, that 
$\varepsilon_0^{(2)}$ is the first non-extraneous root.

The root $\varepsilon_0^{(2)}$ will generally correspond to some 
multiple eigenvalue
$\rho_0$ of $Y_t(\varepsilon_0^{(2)})$, with integers $(p,q)$ as before. Then the statement
of Lemma \ref{yaku} can be formulated as follows. 
\begin{theorem}  \label{thY-S}
If $r \ne \infty$ is the radius of convergence of the series
 \begin{equation}  \label{Y-Seq1}
    \Omega_t(\varepsilon) = \sum_{k=1}^{\infty} \varepsilon^k  \ \Omega_{t,k}, 
 \end{equation}
 there is  at least one non-extraneous root $\varepsilon_0$ of the equation 
 $\Delta(\varepsilon) = 0$ on the circle
 $|\varepsilon| = r$. If for this root one has $p < q$ for some corresponding
 eigenvalue $\rho_0$ of multiplicity $l>1$, then $r = |\varepsilon_0|$, i.e.,  
 the radius of convergence of the series $\Omega_t(\varepsilon)$ is precisely
 $|\varepsilon_0|$.
\end{theorem}
We should remark here that in some cases with $p \ge q$, the series (\ref{Y-Seq1})
may well converge at $\varepsilon = \varepsilon_0$ and the radius of convergence $r$ is
indeed greater than $|\varepsilon_0|$. This occurs, for instance, when $A(t)$ is diagonal. To
illustrate this phenomenon, consider again the matrix (\ref{diag1}) with $a_1(t) \equiv a_2(t)$. 
Then, clearly,
$\rho_1(\varepsilon) = \rho_2(\varepsilon)$ for all $\varepsilon$, so that $l=2$ and $q=1$. If we
choose $\log \rho_1(\varepsilon) = \log \rho_2(\varepsilon)$, then $p=2 > q$.

Although these cases are in a certain sense exceptional, as explained in 
 \cite[page 66]{yakubovich75lde}, Theorem \ref{thY-S}
 is not yet, strictly speaking, a necessary condition for the convergence 
 of the series (\ref{Y-Seq1}). In any
 case, the convergence in the diagonal case is compatible with its formulation, as we
 have seen.

\subsection{Examples}

We next illustrate Theorem \ref{thY-S} on three different examples. We first consider those
analyzed in subsection \ref{examples1} and then we treat in some detail the Magnus
expansion applied to the evolution operator describing a two-level quantum 
system.

\vspace*{0.3cm}

\noindent \textbf{Example 1 (revisited)}. Given the exact solution (\ref{ej1.2}) of Example 1
in subsection \ref{examples1}, the corresponding discriminant is given by
$\Delta(\varepsilon) = (\e^{2 \varepsilon t} + \e^{-\varepsilon t})^2 - 4 \e^{\varepsilon t}$, whose
roots are
\[
     \varepsilon_0^{(1)} = 0, \qquad \mbox{ and }  \qquad \varepsilon_0^{(2)} = i \frac{2 \pi}{3t}.
\]
The first value, $\varepsilon_0^{(1)} = 0$,
 is clearly an extraneous root, so we analyze $\varepsilon_0^{(2)}$.
As $\varepsilon$ varies along the imaginary axis from $\varepsilon = 0$ to 
$\varepsilon = \varepsilon_0^{(2)}$, the eigenvalues of the matrix $Y_t(\varepsilon)$,
\[
   \rho_1(\varepsilon) = \e^{2 \varepsilon t}, \qquad 
   \rho_2(\varepsilon) = \e^{-\varepsilon t}
\]
move along the unit circle, one clockwise and the other counterclockwise from
   \[
       \rho_{1,2}(0) = 1  \quad \mbox{ to } \quad \rho_{1,2}(\varepsilon_0^{(2)}) = 
          \e^{i 4 \pi/3} = \e^{-i 2 \pi/3} = \e^{i(4 \pi/3 - 2 \pi)}.
   \]
Thus, $\rho_1(\varepsilon_0^{(2)})$ and $\rho_2(\varepsilon_0^{(2)})$ lie on different sheets of
  the Riemann surface of the function $\log z$ and therefore $\varepsilon_0^{(2)}$ is a non-extraneous
  root, with $p=1$. Since $Y_t(\varepsilon_0^{(2)})  \ne \rho I$, we have $q=2$, so that,
  according to Theorem \ref{thY-S}, the radius of convergence of the series (\ref{Y-Seq1}) is
  precisely
  \begin{equation}   \label{cex1}
       r = |\varepsilon_0^{(2)}| = \frac{2 \pi}{3t}.
\end{equation}
To get the actual convergence domain of the corresponding Magnus expansion 
we have to take $\varepsilon = 1$, and so, from (\ref{cex1}), we get $2\pi/(3t) = 1$,
or equivalently $t = 2 \pi/3$, i.e., the result achieved from the analysis of the exact solution
in subsection \ref{examples1}.
       
\vspace*{0.3cm}

\noindent \textbf{Example 2 (revisited)}.  Let us obtain the convergence domain for the
Magnus expansion of the solution to the initial value problem $Y' = \varepsilon A(t) Y$, $Y(0)=I$
when $A(t)$ is the piece-wise continuous matrix defined in Example 2 (subsection
\ref{examples1}). The exact solution for $t \ge 1$ is given by
\[
   Y(t;\varepsilon) = \left(  \begin{array}{lc}
        \e^{ \varepsilon w}   &  \ \  \varepsilon \, \beta  \, \e^{\varepsilon w}    \\
        0  &   \e^{-\varepsilon w} 
           \end{array}   \right),
 \]
where $w \equiv \alpha (t-1)$. Equation (\ref{eq.4.1}) leads in this case to
$\cosh^2(\varepsilon w) - 1 = 0$, with first solutions
\[
    \varepsilon = 0, \qquad    \varepsilon = \pm i \frac{\pi}{w}.
\]    
Again, $\varepsilon = 0$ is an extraneous root, whereas the eigenvalues of the
matrix $Y_t(\varepsilon)$ move along the unit circle, one clockwise and the other counterclockwise
from 
   \[
       \rho_{1,2}(0) = 1  \quad  \mbox{ to } \quad \rho_{1,2}(i \pi/w) = -1
   \]
when $\varepsilon$ varies along the imaginary axis from $\varepsilon = 0$ to
$\varepsilon = i \pi/w$ (the same considerations apply to the case $\varepsilon = -i \pi/w$). Then,
obviously, $p=1$ and $q=2$, so that the radius of convergence of the series
(\ref{Y-Seq1}) is 
\[
    |\varepsilon| = \frac{\pi}{|w|} = \frac{\pi}{|\alpha| (t-1)}.
\]
If we now fix $\varepsilon = 1$, we get the actual $t$-domain of convergence
of the Magnus series (\ref{eq.1.4}) as
\[
    t = 1 + \frac{\pi}{|\alpha|}.
\]        
Observe that, when $t=2$, we get $|\alpha| = \pi$ and the result of subsection \ref{examples1}
is recovered: the Magnus
series converges only for $|\alpha| < \pi$.

\vspace*{0.3cm}

\noindent \textbf{Example 3}. Our final illustration corresponds to the quantum mechanical
treatment of a two-level system in a rotating field. It
is described by the Hamiltonian
\begin{equation}   \label{ej3.1}
   H(t) = \frac{1}{2} \hbar \omega_0 \sigma_z + \beta (\sigma_x \cos \omega t + 
   \sigma_y \sin \omega t),
\end{equation}
where $\sigma_x$, $\sigma_y$, $\sigma_z$ are Pauli matrices, and $\beta$ is a coupling
constant.
In fact, this system constitutes a truncation in state space of  a more general one,
namely an atom or freely rotating molecule in a circularly polarized radiation field
\cite{salzman86com,klarsfeld89apo}.  

It has been previously established that when $t = 2 \pi/\omega$ the Magnus
expansion of the corresponding evolution operator $U(t)$, solution of  the
Schr\"odinger equation
\begin{equation}  \label{ej3.2}
   i \hbar \frac{dU}{dt} = H(t) U, \qquad U(0)=I
\end{equation}
converges for $\omega > \omega_0$ and diverges otherwise  
\cite{feldman84otc,salzman86com,maricq87com}. Several different arguments
have been offered trying to explain this phenomenon \cite{klarsfeld89apo}. 
Here we show that this bound can be
directly provided by Theorem \ref{thY-S}.

The exact time-evolution operator  can be obtained in closed form by transforming into a 
rotating frame. Replacing, as usual, $H$ by $\varepsilon H$ in (\ref{ej3.2}) one has
\begin{equation}  \label{ej3.3}
    U(t) = \exp \left(-\frac{1}{2} i \omega t \sigma_z \right) \,
         \exp \left( -i t \Big( \frac{1}{2}(\varepsilon \omega_0 - \omega) \sigma_z +
          \varepsilon \frac{\beta}{\hbar} \sigma_x    \Big) \right).
\end{equation}
From (\ref{ej3.3}), a lengthy but straightforward calculation allows us to write
the corresponding matrix $Y(t;\varepsilon) \equiv U(t)$ in the form
\begin{equation}  \label{ej3.4}
         Y_t(\varepsilon) = \left(  \begin{array}{lc}
           \e^{-\frac{1}{2} i t \omega} \left( \cos \frac{\tilde{\omega} t}{2} - 
                  i \frac{\delta}{\tilde{\omega}} \sin \frac{\tilde{\omega} t}{2} \right) &  
                   \quad -i  \e^{-\frac{1}{2} i t \omega} \frac{2 \varepsilon \beta}{\tilde{\omega} \hbar}
                   \sin \frac{\tilde{\omega} t}{2}  \\
            -i  \e^{\frac{1}{2} i t \omega} \frac{2 \varepsilon \beta}{\tilde{\omega} \hbar}
                   \sin \frac{\tilde{\omega} t}{2}  &  
         \e^{\frac{1}{2} i t \omega} \left( \cos \frac{\tilde{\omega} t}{2} + 
                  i \frac{\delta}{\tilde{\omega}} \sin \frac{\tilde{\omega} t}{2}   \right)      
           \end{array}  \right)
\end{equation}
with $\delta =  \varepsilon \omega_0 - \omega$ and 
$\tilde{\omega} = (\delta^2 + 4 \beta^2 \varepsilon^2/\hbar^2)^{1/2}$. Denoting 
\begin{eqnarray*}
   \cos \chi & \equiv &   \cos \frac{\omega t}{2}   \cos \frac{\tilde{\omega} t}{2} - \frac{\delta}{\tilde{\omega}}
       \sin \frac{\omega t}{2} \sin \frac{\tilde{\omega} t}{2} \\
       & = & \frac{1}{2} (1 + \frac{\delta}{\tilde{\omega}}) \cos \frac{(\omega + \tilde{\omega})t}{2} +        
         \frac{1}{2} (1 - \frac{\delta}{\tilde{\omega}}) \cos \frac{(\omega - \tilde{\omega})t}{2}, 
\end{eqnarray*}
the eigenvalues of $Y_t(\varepsilon)$ can be expressed as
\begin{equation}  \label{ej3.5}
   \rho_{1,2}(\varepsilon) = \cos \chi \pm \sqrt{\cos^2 \chi - 1},
\end{equation}    
so that $U(t)$ has multiple eigenvalues when $\cos \chi = \pm 1$. This equality
is satisfied by $\varepsilon = 0$, which is clearly an extraneous root. The remaining
roots of equation (\ref{eq.4.1}) are obtained from
\begin{equation}  \label{ej3.6}
  \arccos \left(\cos \frac{\omega t}{2}   \cos \frac{\tilde{\omega} t}{2} - \frac{\delta}{\tilde{\omega}}
       \sin \frac{\omega t}{2} \sin \frac{\tilde{\omega} t}{2} \right) = \pi.
\end{equation}       
To simplify the discussion, let us consider the perturbative approximation
$\beta \ll \hbar |\delta|/2$. Then $\tilde{\omega} \approx |\delta|$ and
$\cos \chi \approx \cos(\varepsilon \omega_0 t/2)$, so that (\ref{ej3.6}) reduces to
$\varepsilon \omega_0 t/2 = \pi$. The solution
\[
  \varepsilon_0 = \frac{2 \pi}{\omega_0 t}
\]
is a non-extraneous root with $p=1$ and $q=2$ ($\rho_{1,2}(\varepsilon_0) = -1$), 
and thus the radius of convergence of the
series (\ref{Y-Seq1}) is precisely $|\varepsilon_0|$. 
 Taking now $\varepsilon = 1$, we get finally
the $t$-domain of convergence of the Magnus expansion $t_c = 2\pi/\omega_0$. Notice
that for $t = t_c$ and $\omega < \omega_0$ we are outside the convergence disk, and thus
the Magnus series diverges, just as noted in \cite{klarsfeld89apo}.

\section{Discussion}

The Magnus expansion was originally designed by requiring only that $A(t)$ be a linear
operator depending on a real variable $t$  in an associative ring and that ``certain
unspecified conditions of convergence be satisfied" \cite{magnus54ote}. The idea was to define,
in terms of $A$, an operator $\Omega(t)$ such that 
the solution of the initial value problem
\[
       \frac{dY}{dt} = A(t) Y, \qquad Y(0)=I,
\]
for a second operator $Y$ is given as $Y = \exp \Omega$.
The proposed expression for $\Omega$ was an infinite series satisfying the condition that
 ``its partial sums become Hermitian after multiplication by $i$ if $i A$ is a Hermitian operator"
\cite{magnus54ote}. The simplest example of an equation of this  type is given by
a finite system of linear differential equations. In this case, $A(t)$ is the 
coefficients matrix of the system, and the existence of $\Omega$ is assured ``for
sufficiently small values of $t$" \cite{magnus54ote}. 
Theorem \ref{con-mag} yields an implicit condition in 
terms of the eigenvalues of the matrix $\Omega$. 

Given the importance of the expansion, it has been rediscovered a number of times in
different settings along the years. 
Also a particular attention has been payed to its convergence in the matrix
case, and so
several bounds on the actual radius of convergence of the form (\ref{eq.2.1}) 
have been obtained with different values of $r_c$.  Recently, it has been shown that
the optimal value is $r_c = \pi$ for all $n \times n$ real matrices $A(t)$ \cite{moan06cot}.

In this work, by applying standard techniques of complex analysis and some elementary
properties of the unit sphere, we have generalized this result to bounded linear operators in a Hilbert space (Theorem \ref{conv-mag}), in the spirit of the original Magnus formulation in 
the context of Quantum Mechanics. Obviously, this theorem is also valid for finite dimensional complex
matrices. In our treatment, a complex parameter $\varepsilon$ is introduced in the formalism,
so that the initial value problem (\ref{eq.3.1}) is considered instead. Notice that the Magnus
expansion is trivially recovered as soon as we fix $\varepsilon = 1$.

Although Theorem \ref{conv-mag} provides the optimal convergence domain, in the sense
that  $\pi$ is the largest constant for which the result holds without any further restrictions on 
the operator $A(t)$, one can easily construct examples showing that $\int_0^T \|A(s)\| ds < \pi$
is \emph{not} necessary for the convergence of the expansion. 

With the aim of obtaining a more precise characterization of the convergence, we have
considered in section 4 the case of $n \times n$ complex matrices. There, as a straightforward
consequence of the theory of analytic matrix functions, and in particular, of the logarithm of an 
analytic matrix function such as is done in \cite{yakubovich75lde}, we have established a
connection between the convergence of the Magnus series and the existence of multiple
eigenvalues of the fundamental matrix $Y(t; \varepsilon)$ for a fixed $t$, denoted by
$Y_t(\varepsilon)$
(Theorem \ref{thY-S}). In essence,  if the analytic matrix function $Y_t(\varepsilon)$
has an eigenvalue $\rho_0(\varepsilon_0)$ of multiplicity $l > 1$ 
for a certain $\varepsilon_0$ such that: (a) there is a curve in the $\varepsilon$-plane
joining $\varepsilon = 0$ with $\varepsilon = \varepsilon_0$, and
(b) the number of equal terms in 
$\log \rho_1(\varepsilon_0)$, $\log \rho_2(\varepsilon_0)$, $\ldots, 
\log \rho_{l}(\varepsilon_0)$ such that  $\rho_k(\varepsilon_0) = \rho_0$, $k=1,\ldots, l$ is less
than the maximum dimension of the elementary Jordan block corresponding to $\rho_0$, then
the radius of convergence of the series $\Omega_t(\varepsilon) = 
\sum_{k \ge 1} \varepsilon^k \Omega_{t,k}$ verifying
$\exp \Omega_t(\varepsilon) = Y_t(\varepsilon)$ is precisely $r = |\varepsilon_0|$.

This value $r$ in general will be different for each particular  $t$
considered, so that we can write $r= |\varepsilon_0| = F(t)$ for a given function $F(t)$. In particular,
for the examples considered in section 4, $F(t) = \frac{2 \pi}{3t}$, $F(t) = \frac{\pi}{|\alpha|(t-1)}$ and
$F(t) = \frac{2 \pi}{\omega_0 t}$, respectively. If we fix $|\varepsilon_0| = 1$, then the 
convergence $t$-domain of the Magnus expansion is obtained as the solution of $F(t)=1$ with the
smallest absolute value.

It is interesting at this point to discuss Theorem \ref{con-mag} on the existence of a differentiable
function $\Omega(t)$ and Theorem \ref{conv-mag} in view of the more precise 
account on the convergence
issue provided by Theorem \ref{thY-S}. First, note that under the assumptions of Theorem
\ref{conv-mag}, all the eigenvalues of $Y(t)$ lie in the region 
\[
   G_{\pi} = \{ z = |z| \e^{i \omega} \in \mathbb{C} : \e^{-\pi} \le |z| \le \e^{\pi}, \  
   |\omega| < \pi \},
\]
so that automatically all the differences between any two of the eigenvalues of $\Omega(t) = 
\log Y(t)$ is less than $2 \pi i$ and thus Theorem \ref{con-mag} holds. Second, if all eigenvalues
of $Y_t(\varepsilon)$ are located in $G_{\pi}$, the (possible) multiple eigenvalues take place
only at extraneous roots of the parameter $\varepsilon$ and thus, according with Theorem     
\ref{thY-S}, the convergence of the series is assured.

One might think that the practical application of Theorem \ref{thY-S}
to realistic problems is doubtful, since it is necessary to
compute in advance the fundamental matrix $Y(t;\varepsilon)$. In this sense, the alternative
(but more conservative)
estimate provided by Theorem \ref{conv-mag} directly in terms of the operator $A(t)$ 
is certainly easier to check in practice. In our opinion, however, the characterization of the
convergence of the Magnus expansion in terms of the multiple eigenvalues of $Y(t;\varepsilon)$
sheds new light on this issue, has a theoretical interest by itself and, in addition, 
provides a rigurous justification for the conjecture formulated in 
\cite{moan06cot} on the basis of the exploration of several examples.

\subsection*{Acknowledgements}

The author is grateful to Prof. W. So for providing him the reference \cite{mityagin90unp} and
to Prof. J.A. Oteo for useful discussions. This work has been partially supported by Ministerio de
Educaci\'on y Ciencia (Spain) under project MTM2004-00535
(co-financed by the ERDF of the European Union) and Fundaci\'o Bancaixa.

\bibliographystyle{plain}
\bibliography{ourbib,geom_int,numerbib}

\end{document}